\documentclass[11pt,reqno]{amsart}
\usepackage{amscd,graphicx,psfrag}
\usepackage{amssymb}
\usepackage{color}

\newcommand{\w}{\widehat}
\newcommand{\X}{X}

\newcommand{\p}{\partial}

\newcommand{\B}{\mathcal B}
\newcommand{\Z}{\mathbb Z}

\newcommand{\ep}{\varepsilon}
\renewcommand{\phi}{\varphi}

\newcommand{\ad}{\operatorname{ad}}
\newcommand{\Hom}{\operatorname{Hom}\,}

\newcommand{\tr}{\operatorname{tr}\,}
\newcommand{\Aut}{\operatorname{Aut}}
\newcommand{\Diff}{\operatorname{Diff}}
\newcommand{\lk}{\operatorname{lk}\,}

\renewcommand{\wp}{\hat p\,}
\renewcommand{\i}{\mathbf i}

\newtheorem{theorem}{Theorem}[section]
\newtheorem{lemma}[theorem]{Lemma}
\newtheorem{proposition}[theorem]{Proposition}
\newtheorem{corollary}[theorem]{Corollary}

\newtheorem*{theorem1}{Theorem 1}
\newtheorem*{theorem2}{Theorem 2}

\theoremstyle{definition}
\newtheorem{remark}[theorem]{Remark}
\newtheorem*{example}{Example}

\title{A Casson-Lin type invariant for links}
\thanks{The second author was partially supported by NSF Grant 0305946 and 
the Max-Planck-Institut f\"ur Mathematik in Bonn, Germany}
\author[Eric Harper]{Eric Harper}
\author[Nikolai Saveliev]{Nikolai Saveliev}
\address{Department of Mathematics\newline\indent
University of Miami \newline\indent PO Box 249085
\newline\indent Coral Gables, FL 33124}
\email{\rm{harper@math.miami.edu}}
\email{\rm{saveliev@math.miami.edu}}


\begin{document}

\begin{abstract}
In 1992, Xiao-Song Lin constructed an invariant $h(K)$ of knots $K \subset 
S^3$ via a signed count of conjugacy classes of irreducible $SU(2)$ 
representations of $\pi_1 (S^3 - K)$ with trace-free meridians. Lin showed 
that $h(K)$ equals one half times the knot signature of $K$. Using methods 
similar to Lin's, we construct an invariant $h(L)$ of two-component links 
$L \subset S^3$. Our invariant is a signed count of conjugacy classes of
projective $SU(2)$ representations of $\pi_1 (S^3 - L)$ with a fixed 
$2$-cocycle and corresponding non-trivial $w_2$. We show that $h(L)$ is, 
up to a sign, the linking number of $L$.
\end{abstract}

\maketitle
 

\section{Introduction}
One of the characteristic features of the fundamental group of a closed 
$3$-manifold is that its representation variety in a compact Lie group 
tends to be finite, in a properly understood sense.  This has been a 
guiding principle for defining invariants of $3$-manifolds ever since 
Casson defined his $\lambda$-invariant for integral homology $3$-spheres 
via a signed count of the $SU(2)$ representations of the fundamental group, 
where signs were determined using Heegaard splittings.

Among numerous generalizations of Casson's construction, we will single 
out the invariant of knots in $S^3$ defined by Xiao-Song Lin \cite{Lin} via 
a signed count of $SU(2)$ representations of the fundamental group of 
the knot exterior. The latter is a $3$-manifold with non-empty 
boundary so the above finiteness principle only applies after one imposes 
a proper boundary condition. Lin's choice of the boundary condition,
namely, that all of the knot meridians are represented by trace-free 
$SU(2)$ matrices, resulted in an invariant $h(K)$ of knots $K\subset S^3$. 
Lin further showed that $h(K)$ equals half the knot signature of $K$.

The signs in Lin's construction were determined using braid representations
for knots. Austin (unpublished) and Heusener and Kroll \cite{HK} 
extended this construction by letting the meridians of the knot be 
represented by $SU(2)$ matrices with a fixed trace which need not be zero. 
Their construction gives, for each choice of the trace, a knot invariant 
which equals one half times the equivariant knot signature.

In this paper, we extend Lin's construction to two-component links $L$ in 
$S^3$. In essence, we replace the count of $SU(2)$ representations with a 
count of \emph{projective} $SU(2)$ representations of $\pi_1(S^3 - L)$, in 
the sense of \cite{RS}, with a fixed 2-cocycle representing a non-trivial 
element in the second group cohomology of $\pi_1 (S^3 - L)$. The resulting 
signed count is denoted by $h(L)$. The two main results of this paper are 
then as follows.

\begin{theorem1}
For any two-component link $L \subset S^3$, the integer $h(L)$ is a well 
defined invariant of $L$.
\end{theorem1}

\begin{theorem2}
For any two-component link $L = \ell_1\cup\ell_2$ in $S^3$, one has
\[
h(L)\; =\; \pm \lk(\ell_1,\ell_2).
\] 
\end{theorem2}

It is worth mentioning that our choice of the 2-cocycle imposes Lin's
trace-free condition on us. This is in contrast to Lin's construction, 
where the choice of boundary condition seemed somewhat arbitrary. This 
also means one should not expect to extend our construction to $SU(2)$ 
representations with non-zero trace boundary condition.

Shortly after Casson introduced his invariant for homology $3$-spheres, 
Taubes \cite{Taubes} gave a gauge theoretic description of it in terms 
of a signed count of flat $SU(2)$ connections. After Lin's work, 
but before Heusener and Kroll, a gauge theoretic 
interpretation of the Lin invariant was given by Herald \cite{Herald}. 
He used this interpretation to define an extension of the Lin 
invariant, now known as the Herald--Lin invariant, to knots in 
arbitrary homology spheres, with arbitrary fixed-trace (possibly non-zero) 
boundary condition.

Another attractive feature of the gauge theoretic approach is that it 
can be used to produce ramified versions of the above invariants. Floer 
\cite{Floer} introduced the instanton homology theory whose Euler 
characteristic is twice the Casson invariant. We expect that our 
invariant will have a similar interpretation, perhaps along the lines of 
the knot instanton homology theory of Kronheimer and Mrowka \cite{KM}, 
which in turn is a variant of the orbifold Floer homology of Collin and 
Steer \cite{Collin-Steer}. We hope to discuss this elsewhere, together 
with possible extensions to links in homology spheres and to links of 
more than two components. 


\section{Braids and representations}
Let $F_n$ be a free group of rank $n \ge 2$, with a fixed generating set 
$x_1,\ldots, x_n$. We will follow the conventions of \cite{Long} and 
define the $n$-string braid group $\B_n$ to be the subgroup of $\Aut 
(F_n)$ generated by the automorphisms $\sigma_1,\ldots,\sigma_{n-1}$, 
where the action of $\sigma_i$ is given by 
\[
\begin{array}{llll} 
\sigma_i: & x_i & \mapsto & x_{i+1} \\
& x_{i +1} & \mapsto & (x_{i+1})^{-1}\, x_i\; x_{i+1} \\
& x_j & \mapsto & x_j, \;\; j \neq i ,i +1.
\end{array}
\]

\medskip\noindent
The natural homomorphism $\B_n \to S_n$ onto the symmetric group on $n$ 
letters, $\sigma \mapsto \bar\sigma$, maps each generator $\sigma_i$ to 
the transposition $\bar{\sigma_i}=(i, i + 1)$. A useful observation is 
that, for any $\sigma \in \Aut(F_n)$, one has 
\begin{equation}\label{E:conj}
\sigma (x_i) = w\, x_{\bar{\sigma}^{-1}(i)}\, w^{-1}
\end{equation}
for some word $w \in F_n$. One can also observe that  $\sigma$ preserves 
the product $x_1 \cdots x_n$, that is, 
\begin{equation}\label{E:prod1}
\sigma(x_1 \cdots x_n) = x_1 \cdots x_n
\end{equation}


\subsection{$SU(2)$ representations}
Consider the Lie group $SU(2)$ of unitary two-by-two matrices with 
determinant one, i.e. complex matrices
\[
\left(
\begin{array}{rc}
u        & v       \\
-\bar{v} & \bar{u} \\
\end{array}
\right)
\]

\smallskip\noindent
such that $u\bar{u} + v\bar{v} = 1$. We will often identify $SU(2)$ with 
the group $Sp(1)$ of unit quaternions via
\[
\left(
\begin{array}{rc}
u        & v       \\
-\bar{v} & \bar{u} \\
\end{array}
\right)
\quad \mapsto\quad u + vj\, \in \mathbb H.
\]

\medskip
Let $R_n = \Hom(F_n,SU(2))$ be the space of $SU(2)$ representations of $F_n$, 
and identify it with $SU(2)^n$ by sending a representation $\alpha: F_n 
\to SU(2)$ to the vector $(\alpha(x_1),\ldots,\alpha(x_n))$ of $SU(2)$ 
matrices. The above representation $\B_n \to \Aut(F_n)$ then gives rise to 
the representation
\begin{equation}\label{E:rho1}
\rho:\; \B_n \longrightarrow \Diff (R_n)
\end{equation}
via $\rho(\sigma)(\alpha) = \alpha\circ\sigma^{-1}$. We will abbreviate 
$\rho(\sigma)$ to $\sigma$. We will also denote $\X = (X_1,\ldots,X_n) 
\in R_n$ and write $\sigma(\X) = (\sigma(\X)_1,\ldots,\sigma(\X)_n)$.

\begin{example}
For any $(X_1,\ldots,X_n) \in R_n$, we have $\sigma_1(X_1,X_2,X_3,\ldots,X_n) 
= (X_1 X_2 X_1^{-1},X_1,X_3,\ldots,X_n)$.
\end{example}


\subsection{Extension to the wreath product $\Z_2 \wr \B_n$}
The wreath product $\Z_2 \wr \B_n$ is the semidirect product of $\B_n$ with 
$(\Z_2)^n$, where $\B_n$ acts on $(\Z_2)^n$ by permuting the coordinates, 
$\sigma(\ep_1,\ldots,\ep_n) = (\ep_{\bar\sigma (1)},\ldots,\ep_{\bar\sigma(n)})$. 
Thus the elements of $\Z_2 \wr \B_n$ are the pairs $(\ep,\sigma) \in (\Z_2)^n 
\times \B_n$, with the group multiplication law
\[
(\ep,\sigma)\cdot(\ep',\sigma')\;=\;(\ep\sigma(\ep'),\,\sigma\sigma').
\]
The representation (\ref{E:rho1}) can be extended to a representation
\begin{equation}\label{E:rho3}
\rho:\; \Z_2\wr \B_n \longrightarrow \Diff (R_n)
\end{equation}
by defining 
\[
\rho(\ep,\sigma) (\X) = \ep\cdot\sigma(\X) = (\ep_1 \sigma(\X)_1,
\ldots,\ep_n \sigma(\X)_n),
\]
where $\ep_i$ are viewed as elements of the center $\Z_2 = \{\pm 1\}$ of 
$SU(2)$. That (\ref{E:rho3}) is a representation follows by a direct 
calculation after one observes that, because of (\ref{E:conj}),
\begin{equation}\label{E:conj2}
\sigma(\X)_i = A X_{\bar\sigma(i)} A^{-1}\quad\text{for some $A \in SU(2)$}.
\end{equation}
Again, we will abuse notations and write simply $\ep\sigma$ for both 
$(\ep,\sigma)$ and $\rho(\ep,\sigma)$.

\begin{example} For any $(X_1,\ldots,X_n) \in R_n$ and $\ep = (\ep_1,
\ldots,\ep_n) \in (\Z_2)^n$, we have 
$(\ep\sigma_1)(X_1,X_2,X_3,\ldots,X_n) = (\ep_1\,X_1 X_2 
X_1^{-1},\ep_2\,X_1,\ep_3\,X_3,\ldots,\ep_n\,X_n)$.
Also $\sigma_1(\ep X) = \sigma_1(\ep) \sigma_1(\X) 
= (\ep_2 X_1 X_2 X_1^{-1} , \ep_1 X_1, \ep_3 X_3, \ldots , \ep_n X_n)$.   
\end{example}


\subsection{Braids and link groups}
The closure $\hat\sigma$ of a braid $\sigma \in \B_n$ is a link in $S^3$ 
with link group
\[
\pi_1 (S^3 - \hat\sigma) = \langle\,x_1,\ldots, x_n\;|\;x_i = \sigma (x_i),
\; i = 1,\ldots, n\,\rangle,
\]
where each $x_i$ represents a meridian of $\hat\sigma$. One can easily see 
that the fixed points of the diffeomorphism $\sigma: R_n \to R_n$ are 
representations $\pi_1 (S^3 - \hat\sigma) \to SU(2)$. This paper grew out 
of the observation that a fixed point $\alpha = (\alpha(x_1),\ldots,
\alpha(x_n))$ of the diffeomorphism $\ep\sigma: R_n\to R_n$ gives rise to 
a representation $\ad\alpha: \pi_1(S^3 - \hat\sigma) \to SO(3)$ by 
composing with the adjoint representation $\ad: SU(2) \to SO(3)$. 
Depending on $\ep$, the representation $\ad\alpha$ may or may not
lift to an $SU(2)$ representation, the obstruction being the second 
Stiefel--Whitney class $w_2 (\ad\alpha) \in H^2(\pi_1(S^3 - \hat\sigma);
\Z_2)$.


\section{Definition of $h(\ep\sigma)$}
Every link in $S^3$ is the closure $\hat\sigma$ of a braid $\sigma$; see 
Alexander \cite{A}. Let $\sigma$ be a braid whose closure $\hat\sigma$ 
has two components. We will associate with it, for a carefully chosen $\ep$, 
an integer $h(\ep\sigma)$. We will prove in Section \ref{S:link-inv} that 
$h$ is an invariant of the link $\hat\sigma$.


\subsection{Choice of $\ep$}
The number of components of the link $\hat\sigma$ is exactly the number 
of cycles in the permutation $\bar\sigma$. We will be interested in two 
component links, that is, the closures of braids $\sigma$ with 
\begin{equation}\label{E:cycles}
\bar\sigma = (i_1\ldots i_m)(i_{m+1}\ldots i_n)\quad\text{for some\; 
$1\le m\le n - 1$}.
\end{equation}
Given such a braid $\sigma$, choose a vector $\ep = (\ep_1,\ldots,
\ep_n) \in (\Z_2)^n$ such that
\begin{equation}\label{E:ep}
\ep_{i_1} \cdots \ep_{i_m} = \ep_{i_{m+1}} \cdots \ep_{i_n} = -1.
\end{equation}
This choice of $\ep$ is dictated by the following two considerations. 
First, we wish to preserve condition (\ref{E:prod1}) in the form
\begin{equation}\label{E:prod2}
(\ep\sigma)(\X)_1\cdots (\ep\sigma)(\X)_n = X_1\cdots X_n,
\end{equation}
and second, we want the fixed points $\alpha$ of the diffeomorphism 
$\ep\sigma: R_n \to R_n$ to have non-zero $w_2\,(\ad\alpha)$.

\begin{lemma}\label{E:w2}
Let $\alpha$ be a fixed point of $\ep\sigma: R_n\to R_n$ with $\ep$ 
as in (\ref{E:ep}) then $w_2\,(\ad\alpha) \ne 0$. 
\end{lemma}

\begin{proof}
The class $w_2\,(\ad\alpha)$ is the obstruction to lifting $\ad\alpha$ 
to an $SU(2)$ representation. Extend $\alpha$ arbitrarily to a function 
$\alpha: \pi_1 (S^3 - \hat\sigma) \to SU(2)$ lifting $\ad\alpha$ then 
$w_2\,(\ad\alpha)$ will vanish if and only if there is a function $\eta: 
\pi_1 (S^3 - \hat\sigma) \to \Z_2 = \{\pm 1\}$ such that $\eta \cdot 
\alpha$ is a representation. Suppose that such a function exists, and 
denote $\eta(x_i) = \eta_i = \pm 1$. Also, assume without loss of 
generality that $\bar\sigma = (1\ldots m)(m+1 \ldots n)$. It follows from 
(\ref{E:conj2}) that in order to satisfy the relations $X_i = (\ep\sigma)
(\X)_i$ we must have $\eta_1 = \ep_1\eta_2 = \ep_1\ep_2\eta_3 = \ldots = 
\ep_1\cdots\ep_m\eta_1 = - \eta_1$, a contradiction with $\eta_1 = \pm 1$.
\end{proof}

The above result concerning $w_2\,(\ad\alpha)$ can be sharpened using the 
following algebraic topology lemma.

\begin{lemma}
Let $\w\sigma$ be a link of two components. If $\w\sigma$ is non-split then 
$H^2 (\pi_1 (S^3 - \hat\sigma);\Z_2) = \Z_2$. Otherwise, $H^2 (\pi_1 (S^3 - 
\hat\sigma);\Z_2) = 0$.
\end{lemma}

\begin{proof} 
If $\w\sigma$ is non-split then $S^3 - \w\sigma$ is a $K(\pi,1)$ by the 
Sphere Theorem, hence $H^2(\pi_1(S^3 - \w\sigma);\Z_2) = H^2 (S^3 - \w\sigma;
\Z_2) = \Z_2$. If $\w\sigma$ is split then $K(\pi_1(S^3-\w\sigma),1)$ has 
the homotopy type of a one-point union of two circles and the result again 
follows.
\end{proof}

\begin{corollary}\label{C:w2}
Let $\w\sigma$ be a split link of two components, and let $\ep$ be chosen 
as in (\ref{E:ep}). Then the diffeomorphism $\ep\sigma: R_n \to R_n$ has no 
fixed points. 
\end{corollary}


\subsection{The zero-trace condition}
A naive way to define $h(\ep\sigma)$ would be as the intersection number 
of the graph of $\ep\sigma: R_n \to R_n$ with the diagonal in the product 
$R_n \times R_n$. One can observe though that, in addition to this 
intersection not being transversal, its points $(\X,\X) = (X_1,\ldots,X_n,
X_1,\ldots,X_n)$ have the property that $\tr X_1 = \ldots = \tr X_n = 0$. 
This can be seen as follows. 

Assume without loss of generality that $\bar{\sigma} = (1 \ldots m) (m+1 
\ldots n)$. Then the relations $\X = \ep \sigma (\X)$ together with 
(\ref{E:conj2}) imply that 
\begin{multline}\notag
X_1 = \ep_1 \sigma(\X)_1 = \ep_1 A_1\cdot X_{\bar\sigma(1)}\cdot A_1^{-1} 
= \ep_1 A_1 X_2 A_1^{-1} \\ = \ep_1 A_1\cdot \ep_2\,\sigma(\X)_2\cdot A_1^{-1}
= \ep_1\ep_2\,A_1 A_2\cdot X_{\bar\sigma(2)}\cdot A_2^{-1} A_1^{-1}= \ldots \\ 
= \ep_1\cdots\ep_m\,(A_1\cdots A_m)\cdot X_1\cdot (A_1 \ldots A_m)^{-1}.
\end{multline}
Since trace is conjugation invariant and $\ep_1 \ldots \ep_m = -1$, we 
conclude that $\tr X_1 = \ldots = \tr X_m = 0$. Similarly, $\tr X_{m+1} =
\ldots = \tr X_n = 0$.  

Therefore, in our definition we will restrict ourselves to the subset of 
$R_n$ consisting of $\X = (X_1,\ldots,X_n)$ with $\tr X_1 = \ldots = \tr X_n 
= 0$. The non-transversality problem will be addressed below by factoring 
out the conjugation symmetry and lowering the dimension of the ambient 
manifold.


\subsection{The definition}
The subset of $SU(2)$ consisting of the matrices with zero trace is a 
conjugacy class in $SU(2)$ diffeomorphic to $S^2$. Define
\[
Q_n = \{(X_1,\ldots,X_n) \in R_n\; |\; \tr X_i = 0\,\} \subset R_n,
\]
so that $Q_n$ is a manifold diffeomorphic to $(S^2)^n$. Also define 
\[
H_n  = \{(X_1,\ldots,X_n,Y_1,\ldots,Y_n) \in Q_n \times Q_n\; |\; 
X_1\cdots X_n = Y_1\cdots Y_n\,\}.
\]
This is no longer a manifold due to the presence of {\em reducibles}. We 
call a point $(X_1,\ldots,X_n,Y_1,\ldots,Y_n)\in Q_n\times Q_n$ reducible
if all $X_i$ and $Y_j$ commute with each other or, equivalently, if there 
is a matrix $A \in SU(2)$ such that $A X_i A^{-1}$ and $A Y_i A^{-1}$ are 
all diagonal matrices, $i = 1,\ldots,n$.  The subset $S_n \subset Q_n 
\times Q_n$ of reducibles is closed. 

\begin{lemma}
$H_n^* = H_n - S_n$ is an open manifold of dimension $4n-3$.
\end{lemma}

\begin{proof}
Let us consider the open manifold $(Q_n\times Q_n)^* = Q_n\times Q_n - S_n$ 
of dimension $4n$ and the map $f: (Q_n \times Q_n)^* \to SU(2)$ given by
\begin{equation}\label{E:mapf}
f(X_1, \ldots, X_n,Y_1, \ldots, Y_n) = X_1 \cdots X_n Y_n^{-1} \cdots Y_1^{-1}.
\end{equation}
According to Lemma 1.5 of \cite{Lin}, this map has $1 \in SU(2)$ as a regular
value. Since $H_n^* = f^{-1} (1)$, the result follows. 
\end{proof}

Because of (\ref{E:conj2}) and the fact that multiplication by $-1\in SU(2)$ 
preserves the zero trace condition, the representation (\ref{E:rho3}) gives 
rise to a representation
\begin{equation}\label{E:rho4}
\rho:\; \Z_2\wr\B_n \longrightarrow \Diff (Q_n).
\end{equation}
Given $\ep\sigma \in \Z_2 \wr \B_n$ such that (\ref{E:cycles}) and 
(\ref{E:ep}) are satisfied, consider two sub\-manifolds of $Q_n \times Q_n$: 
one is the graph of $\ep\sigma: Q_n \to Q_n$, 
\[
\Gamma_{\ep\sigma} = \{\,(\X,\ep\sigma(\X))\;|\;\X \in Q_n\,\},
\]
and the other the diagonal, 
\[
\Delta_n = \{\,(\X,\X)\;|\;\X \in Q_n\,\}.
\]
Note that both $\Gamma_{\ep\sigma}$ and $\Delta_n$ are subsets of $H_n$\,:
this is obvious for $\Delta_n$, and follows from equation (\ref{E:prod2}) 
for $\Gamma_{\ep\sigma}$.

\begin{proposition}\label{P:irrint}
The intersection $\Gamma_{\ep\sigma}\,\cap\,\Delta_n \subset H_n$ consists 
of irreducible representations.
\end{proposition}   

\begin{proof}
Assume without loss of generality that $\bar\sigma = (1\ldots m)(m+1
\ldots n)$, and suppose that $(\X,\X) = (X_1,\ldots,X_n,X_1,\ldots,X_n)\in
\Gamma_{\ep\sigma}\,\cap\,\Delta_n$ is reducible. Then all  of the $X_i$ commute 
with each other and, in particular, $\sigma(\X) = (X_{\bar\sigma(1)},\ldots,
X_{\bar\sigma(n)})$. The equality $\X = \ep\sigma(\X)$ then implies that 
$X_1 = \ep_1 X_{\bar\sigma(1)} = \ep_1 X_2 = \ep_1\ep_2 X_{\bar\sigma(2)} = 
\ldots = \ep_1\cdots\ep_m X_1 = - X_1$, a contradiction with $X_1\in SU(2)$.
\end{proof}

Let $\Gamma^*_{\ep\sigma} = \Gamma_{\ep\sigma}\,\cap\,H_n^*$ and $\Delta_n^* 
= \Delta_n\,\cap\, H_n^*$ be the irreducible parts of $\Gamma_{\ep\sigma}$
and $\Delta_n$, respectively. They are both open submanifolds of $H_n^*$ of
dimension $2n$.

\begin{corollary}
The intersection $\Delta_n^*\,\cap\,\Gamma_{\ep\sigma}^* \subset H_n^*$ is 
compact.
\end{corollary}

\begin{proof}
Proposition \ref{P:irrint} implies that $\Delta_n^*\,\cap\,
\Gamma_{\ep\sigma}^* = \Delta_n\,\cap\,\Gamma_{\ep\sigma}$, and the latter 
intersection is obviously compact as it is the intersection of two compact 
subsets of $H_n$. 
\end{proof}

The group $SO(3) = SU(2)/\{\pm 1\}$ acts freely by conjugation on $H_n^*$, 
$\Delta_n^*$, and $\Gamma_{\ep\sigma}^*$. Denote the resulting quotient 
manifolds by 
\[
\w H_n = H_n^*/SO(3),\quad \w \Delta_n = \Delta^*_n/SO(3), \quad\text{and}
\quad \w \Gamma_{\ep\sigma} = \Gamma^*_{\ep\sigma}/SO(3). 
\]
The dimension of $\w H_n$ is $4n - 6$, and $\w \Delta_n$ and 
$\w \Gamma_{\ep\sigma}$ are submanifolds, each of dimension $2n-3$.
Since the intersection $\w \Delta_n\,\cap\,\w \Gamma_{\ep\sigma}$ is 
compact, one can isotope $\w \Gamma_{\ep\sigma}$ into a submanifold 
$\widetilde \Gamma_{\ep\sigma}$ using an isotopy with compact support 
so that $\w \Delta_n\,\cap\,\widetilde \Gamma_{\ep\sigma}$ consists of 
finitely many points. Define 
\[
h(\ep\sigma) = \#_{\w H_n}\,(\w\Delta_n\,\cap\,\widetilde\Gamma_{\ep\sigma})  
\]
as the algebraic intersection number, where the orientations of $\w H_n$, 
$\w \Delta_n$, and $\widetilde\Gamma_{\ep\sigma}$ are described in the 
following subsection. It is obvious that $h(\ep\sigma)$ does not depend on 
the perturbation of $\w \Gamma_{\ep\sigma}$ so we will simply write
\[
h(\ep\sigma) = \langle\,\w \Delta_n, \w \Gamma_{\ep\sigma}\,\rangle_{\w H_n}.
\]


\subsection{Orientations}
Orient the copy of $S^2 \subset SU(2)$ cut out by the trace zero condition 
arbitrarily, and endow $Q_n = (S^2)^n$ and $Q_n \times Q_n$ with product 
orientations. The diagonal $\Delta_n$ and the graph $\Gamma_{\ep\sigma}$ 
are naturally diffeomorphic to $Q_n$ via projection onto the first factor, 
and they are given the induced orientations. Note that if we reverse the 
orientation of $S^2$, then the orientation of $Q_n$ is reversed if $n$ is 
odd. Hence the orientations of both $\Delta_n$ and $\Gamma_{\ep\sigma}$ are 
reversed if $n$ is odd, while the orientation of $Q_n \times Q_n = 
(S^2)^{2n}$ is preserved regardless of the parity of $n$. 

Orient $SU(2)$ by the standard basis $\{ i, j, k \}$ in its Lie algebra 
$\mathfrak{su}(2)$, and orient $H_n^* = f^{-1}(1)$ by applying the 
base--fiber rule to the map (\ref{E:mapf}). The adjoint action of $SO(3)$ 
on $S^2 \subset SU(2)$ is orientation preserving, hence the $SO(3)$ 
quotients $\w H_n$, $\w \Delta_n$, and $\w \Gamma_{\ep\sigma}$ are 
orientable. We orient them using the base--fiber rule. The discussion in 
the previous paragraph shows that reversing orientation on $S^2$ may 
reverse the orientations of $\w \Delta_n$ and $\w \Gamma_{\ep\sigma}$ but 
that it does not affect the intersection number $\langle\, \w \Delta_n, 
\w \Gamma_{\ep\sigma}\, \rangle_{\w H_n}$.


\section{The link invariant $h$}\label{S:link-inv}
In this section, we will prove Theorem 1. This will be accomplished by 
proving that $h(\ep\sigma)$ is independent first of $\ep$ and then of 
$\sigma$.


\subsection{Independence of $\ep$}
We will first show that, for a fixed $\sigma$ whose closure $\hat\sigma$ 
is a link of two components, $h(\ep\sigma)$ is independent of the choice of 
$\ep$ as long as $\ep$ satisfies (\ref{E:ep}).

\begin{proposition}\label{E:indep}
Let $\ep$ and $\ep'$ be such that (\ref{E:ep}) is satisfied. Then 
$h(\ep\sigma) = h(\ep'\sigma)$.
\end{proposition}

\begin{proof}
Assume without loss of generality that $\bar\sigma = (1\ldots m)
( m+1 \ldots n)$ and let $\ep = (\ep_1, \ldots, \ep_n)$ and $\ep' =
(\ep'_1, \ldots, \ep'_n)$. Define $\delta = (\delta_1, \ldots, \delta_n)$ 
as the vector in $(\Z_2)^n$ with coordinates
\[
\delta_1 = 1\quad\text{and}\quad \delta_{k+1} \;=\; \delta_k\,\ep_k\,\ep'_k 
\quad\text{for}\quad k = 1,\ldots, n-1,
\]
and define the involution $\tau: Q_n \to Q_n$ by the formula
\[
\tau (X) = \delta X = (\delta_1 X_1,\delta_2 X_2,\ldots, 
\delta_n X_n).
\]
Recall that $Q_n = (S^2)^n$ so that $\tau$ is a diffeomorphism which 
restricts to each of the factors $S^2$ as either the identity or the 
antipodal map. In particular, $\tau$ need not be orientation preserving.

The map $\tau \times \tau: Q_n \times Q_n \rightarrow Q_n \times Q_n$
obviously preserves the irreducibility condition and commutes with the 
$SO(3)$ action. It gives rise to an orientation preserving automorphism 
of $\w H_n$ which will again be called $\tau\times\tau$. It is clear 
that $(\tau \times \tau)(\w \Delta_n) = \w \Delta_n$. It is also true 
that $(\tau \times \tau)(\w \Gamma_{\ep\sigma}) = \w\Gamma_{\ep'\sigma}$, 
which can be seen as follows. Given a pair $(\delta X,\delta\ep\sigma(X))$ 
whose conjugacy class belongs to $(\tau\times\tau)(\w\Gamma_{\ep\sigma})$,
write it as
\[
(\delta X, \delta \ep \sigma (X)) =
(\delta X, \delta \ep \sigma (\delta \delta X)) = 
(\delta X, \delta \ep \sigma(\delta)\,\sigma(\delta X)) 
\] 
using the multiplication law in the group $\Z_2\wr\B_n$. The conjugacy
class of this pair belongs to $\Gamma_{\ep'\sigma}$ if and only 
if $\delta \ep \sigma(\delta) = \ep'$. That this condition holds can be 
verified directly from the definition of $\delta$.

Recall that the orientations of $\w \Delta_n$, $\w \Gamma_{\ep \sigma}$, and 
$\w \Gamma_{\ep' \sigma}$ are induced by the orientation of $Q_n$. Therefore, 
the maps $\tau \times \tau: \w \Delta_n \to \w \Delta_n$ and $\tau \times 
\tau: \w \Gamma_{\ep \sigma} \to \w \Gamma_{\ep' \sigma}$ are either both 
orientation preserving or both orientation reversing depending on 
whether $\tau: Q_n \to Q_n$ preserves or reverses orientation. Hence
we have
\begin{multline}\notag
h(\ep\sigma)
=  \langle\,\w \Delta_n , \w \Gamma_{\ep\sigma}\,\rangle_{\w H_n}
=  \langle\,(\tau \times \tau) (\w \Delta_n), (\tau \times \tau) 
(\w \Gamma_{\ep\sigma})\,\rangle_{(\tau \times \tau)(\w H_n)} \\
=  \langle\,\w \Delta_n , \w \Gamma_{\ep'\sigma}\,\rangle_{\w H_n}
=  h(\ep'\sigma).   
\end{multline}
\end{proof}

From now on, we will drop $\ep$ from the notation and simply write 
$h(\sigma)$ for $h(\ep\sigma)$ assuming that a choice of $\ep$ 
satisfying (\ref{E:ep}) has been made. 


\subsection{Independence of $\sigma$}
In this section, we will show that $h(\sigma)$ only depends on the link 
$\hat \sigma$, not on a particular choice of braid $\sigma$, by 
verifying that $h$ is preserved under Markov moves. We will follow the 
proof of \cite[Theorem 1.8]{Lin} which goes through with little change 
once the right $\ep$ are chosen. 

Recall that two braids $\alpha \in \B_n$ and $\beta$ $\in$ $\B_m$ have 
isotopic closures $\hat \alpha$ and $\hat \beta$ if and only if one braid 
can be obtained from the other by a finite sequence of Markov moves; see 
for instance \cite{Birman}. A type 1 Markov move replaces $\sigma \in 
\B_n$ by $\xi^{-1} \sigma \xi \in \B_n$ for any $\xi \in \B_n$. A type 2 
Markov move means replacing $\sigma \in \B_n$ by $\sigma_n^{\pm 1}\sigma 
\in \B_{n+1}$, or the inverse of this operation.

\begin{proposition}
The invariant $h(\sigma)$ is preserved by type 1 Markov moves.
\end{proposition}

\begin{proof}
Let $\xi$, $\sigma \in \B_n$ and assume as usual that $\bar{\sigma} = 
(1 \ldots m)(m+1 \ldots n)$. Then 
\[
\overline{\xi^{-1}\sigma\xi} = 
(\bar{\xi}(1) \ldots \bar{\xi}(m)) (\bar{\xi}(m+1) \ldots \bar{\xi}(n))
\]
has the same cycle structure as $\bar\sigma$. To compute $h(\xi^{-1}\sigma
\xi)$, we will make a choice of $\ep \in (\Z_2)^n$ which satisfies
condition (\ref{E:ep}) with respect to the braid $\xi^{-1}\sigma\xi$, that 
is, $\ep_{\bar\xi(1)}\cdots \ep_{\bar\xi(m)} = \ep_{\bar\xi(m+1)}\cdots 
\ep_{\bar\xi(n)} = -1$.

The braid $\xi$ gives rise to the map $\xi: Q_n \rightarrow Q_n$. It acts 
by permutation and conjugation on the $S^2$ factors in $Q_n$ hence it is 
orientation preserving (we use the fact that $S^2$ is even dimensional). 
It induces an orientation preserving map $\xi \times \xi: Q_n \times Q_n 
\rightarrow Q_n \times Q_n$, which preserves the irreducibility condition 
and commutes with the $SO(3)$ action. Equation (\ref{E:prod1}) then 
ensures that we have a well defined orientation preserving automorphism 
$\xi \times \xi: \w H_n \to \w H_n$. 

That this automorphism preserves the diagonal, $(\xi \times \xi) 
(\w \Delta_n) = \w \Delta_n$, is obvious. Concerning the graphs, let 
$(\X,\ep\xi^{-1}\sigma\xi(\X)) \in \w \Gamma_{\ep\xi^{-1}\sigma\xi}$ then 
\begin{multline}\notag
(\xi \times \xi) (\X, \ep\xi^{-1}\sigma\xi(\X)) \\ = 
(\xi(\X),\xi(\ep\xi^{-1}\sigma\xi(\X))) =
(\xi(\X),\xi(\ep)\sigma(\xi(\X))) \in \w \Gamma_{\xi(\ep)\sigma}.
\end{multline}
Therefore, $(\xi \times \xi)(\w \Gamma_{\ep\xi^{-1}\sigma\xi}) = 
\w \Gamma_{\xi(\ep)\sigma}$. Since $\xi: Q_n \to Q_n$ is orientation 
preserving, the above identifications of the diagonals and graphs via $\xi
\times \xi$ are also orientation preserving. 

Observe that $\xi(\ep)_i = \ep_{\bar\xi(i)}$ hence $\xi(\ep)$ satisfies 
(\ref{E:ep}) with respect to $\sigma$ and thus can be used to compute 
$h(\sigma)$. The following calculation now completes the argument\,:
\begin{multline}\notag
h(\xi^{-1}\sigma\xi) = 
\langle \w \Delta_n , \w \Gamma_{\ep\xi^{-1}\sigma\xi} \rangle_{\w H_n}
=  \langle (\xi \times \xi) (\w \Delta_n) , (\xi \times \xi) 
(\w \Gamma_{\ep\xi^{-1}\sigma\xi}) \rangle_{(\xi \times \xi) (\w H_n)} \\
=  \langle \w \Delta_n , \w \Gamma_{\xi(\ep)\sigma} \rangle_{\w H_n} 
=  h(\sigma).
\end{multline}
\end{proof}

\begin{proposition}
The invariant $h(\sigma)$ is preserved by type 2 Markov moves.
\end{proposition}

\begin{proof}
Given $\sigma \in \B_n$ and $\ep$ satisfying (\ref{E:ep}), change $\sigma$ 
to $\sigma_n\sigma \in \B_{n+1}$ and let $\ep' = \sigma_n(\ep,1)$.  If $X = 
(X_1,\ldots,X_n)$ and $\ep = (\ep_1,\ldots,\ep_n)$ then 
\begin{multline}\notag
(\sigma_n\sigma)(X,X_{n+1}) = \sigma_n (\sigma(X),X_{n+1}) \\ 
= (\sigma(X)_1,\ldots,\sigma(X)_{n-1},\,\sigma(X)_n X_{n+1}\sigma(X)_n^{-1},
\,\sigma(X)_n))
\end{multline}
and $\ep' = (\ep_1,\ldots,\ep_{n-1},1,\ep_n)$. In particular, $\ep'$ satisfies
(\ref{E:ep}) with respect to $\sigma_n\sigma$. Consider the embedding $g:
Q_n \times Q_n \rightarrow Q_{n+1} \times Q_{n+1}$ given by
\[
g(X_1, \ldots, X_n, Y_1, \ldots, Y_n) 
= (X_1, \ldots, X_n, Y_n, Y_1, \ldots, Y_n, Y_n)
\]
One can easily see that $g (H_n) \subset H_{n+1}$ and that $g$ commutes with 
the conjugation, thus giving rise to an embedding $\hat g: \w H_n \to 
\w H_{n+1}$. A straightforward calculation using the above formulas for 
$\sigma_n\sigma$ and $\ep'$ then shows that 
\[
\quad \hat g\,(\w\Delta_n) \subset \w\Delta_{n+1},\quad
\hat g\,(\w\Gamma_{\ep\sigma}) \subset \w\Gamma_{\ep'\sigma_n\sigma}\, ,\;\;
\text{and}\quad \hat g\,(\w\Delta_n\,\cap\, \w \Gamma_{\ep\sigma}) = 
\w\Delta_{n+1}\,\cap\, \w \Gamma_{\ep'\sigma_n\sigma}.
\]
Now, one can achieve all the necessary transversalities and match the 
orientations in exactly the same way as in the second half of the proof 
of \cite[Theorem 1.8]{Lin}. This shows that $h(\sigma_n\sigma) = h(\sigma)$. 
The proof of the equality $h(\sigma_n^{-1}\sigma) = h(\sigma)$ is similar.
\end{proof}


\section{The invariant $h(\sigma)$ as the linking number}
In this section we will prove Theorem 2, that is, show that for any link 
$\hat\sigma = \ell_1 \cup \ell_2$ of two components, one has 
\[
h(\sigma)\; =\; \pm\, \lk (\ell_1,\ell_2).
\]
Our strategy will be to show that the invariant $h(\sigma)$ and the 
linking number $\lk (\ell_1,\ell_2)$ change according to the same rule as 
we change a crossing between two strands from two different components of 
$\hat\sigma = \ell_1 \cup \ell_2$ (the link $\hat\sigma$ will need to be 
oriented for that, although a particular choice of orientation will not 
matter). After changing finitely many such 
crossings, we will arrive at a split link, for which both the invariant 
$h(\sigma)$ and the linking number $\lk (\ell_1,\ell_2)$ vanish; see 
Corollary \ref{C:w2}. The change of crossing as above obviously changes 
the linking number by $\pm 1$. To calculate the effect of the crossing 
change on $h(\sigma)$, we will follow \cite{Lin} and reduce the problem to 
a calculation in the pillowcase $\w{H}_2$.


\subsection{The pillowcase}
\begin{figure}[!ht]
\centering
\psfrag{psi}{$\psi$}
\psfrag{pi}{$\pi$}
\psfrag{t}{$\theta$}
\psfrag{0}{$0$}
\psfrag{G}{$\w\Gamma_{\ep\sigma_1^2}$}
\psfrag{D}{$\w\Delta_2$}
\includegraphics{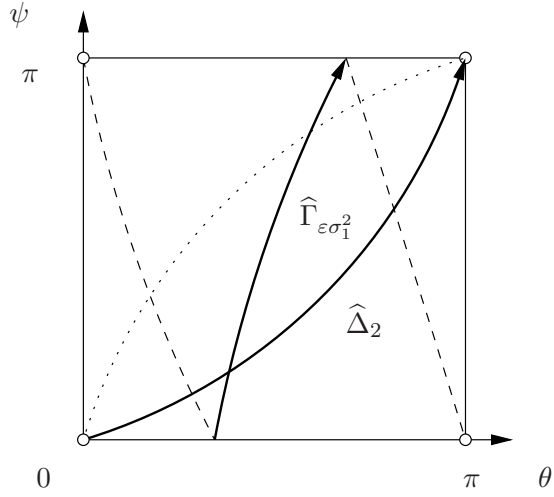}
\caption{The pillowcase}
\label{E:hopflink}
\end{figure}

We begin with a geometric description of $\w H_2$ as a pillowcase, compare
with \cite[Lemma 1.2]{Lin}. Remember that 
\[
H_2 = \{(X_1, X_2,Y_1, Y_2) \in Q_2 \times Q_2\; |\; X_1 X_2 
= Y_1 Y_2\,\}.
\]
We will use the identification of $SU(2)$ with $Sp(1)$ when convenient. 
Since $X_2$ is trace free, we may assume that $X_2 = i$ after conjugation. 
Conjugating by $e^{i\phi}$ will not change $X_2$ but, for an 
appropriate choice of $\phi$, will make $X_1$ into 
\[
X_1=
\begin{pmatrix} ir & u \\ - u & -ir \end{pmatrix},
\]
where both $r$ and $u$ are real, and $u$ is also non-negative. Since $r^2 + 
u^2 = 1$ we can write $r = \cos\theta$ and $u = \sin\theta$ for a unique 
$\theta $ such that $0 \leq \theta \leq \pi$. In the quaternionic language, 
$X_1 = i\, e^{-k \theta}$ with $0 \leq \theta \leq \pi$. Similarly, the 
condition $\tr (Y_2) = \tr(Y_1^{-1} X_1 X_2) = 0$ implies that $Y_1 = 
i\, e^{-k \psi}$, this time with $-\pi \leq \psi \leq \pi$. To summarize,
\[
X_1 = i e ^{-k \theta},\;\; X_2 = i,\;\; Y_1 =  ie ^{-k \psi},\;\; 
Y_2 = i e^{-k(\psi -\theta)}.
\]
Thus $\w H_2$ is parameterized by the rectangle $[0,\pi]
\times [-\pi,\pi]$, with proper identifications along the edges and with 
the reducibles removed. The reducibles occur when both $\theta$ and $\psi$ 
are multiples of $\pi$, hence $\w H_2$ a 2-sphere with the points $A 
= (0,0)$, $B = (\pi,0)$, $A' = (0,\pi)$, and $B' = (\pi,\pi)$ removed; see 
Figure \ref{E:hopflink}. According to \cite{Lin}, the orientation on the 
front sheet of $\w H_2$ coincides with the standard orientation on the 
$(\theta,\psi)$ plane. 

\begin{example}\label{E:hopf}
Let $\sigma = \sigma_1^2$ so that $\w \sigma = \ell_1 \cup \ell_2$ is the 
Hopf link with $\lk (\ell_1,\ell_2) = \pm 1$. To calculate $h(\sigma)$, we 
let $\ep = (-1,-1)$, the only available choice satisfying (\ref{E:ep}), and 
consider the submanifolds $\w{\Delta}_2$ and $\w{\Gamma}_{\ep\sigma}$ of
$\w H_2$. We have, in quaternionic notations, $\w{\Delta}_2 = 
\{ (i e^{-k \theta}, i, i e ^{-k \theta}, i)\}$, which is the diagonal 
$\psi = \theta$ in the pillowcase. A straightforward calculation shows that 
$\w {\Gamma}_{\ep\sigma} \subset \w H_2$ is given by $\psi = 3\theta - \pi$.  
As can be seen in Figure \ref{E:hopflink}, the intersection $\w \Delta_2 
\,\cap\,\w\Gamma_{\ep\sigma}$ consists of one point coming with a sign. 
Hence $h(\sigma_1^2) = \pm 1$, which is consistent with the fact that 
$\lk (\ell_1,\ell_2) = \pm 1$.
\end{example}

\begin{example}
Let $\sigma = \sigma_1^{2n}$ then arguing as above one can show that
$\w \Gamma_{\ep\sigma}\subset \w H_2$ is given by $\psi = (2n+1)\theta - 
\pi$. In this case there are $n$ intersection points all of which come 
with the same sign. This shows that $h(\sigma_1^{2n}) = \pm n$, which is 
again consistent with the fact that $\lk (\ell_1,\ell_2) = \pm n$.
\end{example}


\subsection{The difference cycle}
Given a two component link $\hat\sigma$, fix an orientation on it. A 
particular choice of orientation will not matter because we are only
interested in identifying $h(\sigma)$ with the linking number $\lk 
(\ell_1,\ell_2)$ up to sign. We wish to change one of the crossings 
between the two components of $\hat\sigma$. Using a sequence of first 
Markov moves, we may assume that the first two strands of $\sigma$ 
belong to two different components of $\hat\sigma$, and that the 
crossing change occurs between these two strands. Furthermore, we may 
assume that the crossing change makes $\sigma$ into $\sigma_1^{\pm 2}
\sigma$, where the sign depends on the type of the crossing we change. 
Note that the braid $\sigma_1^{\pm 2}\sigma$ has the same permutation 
type as $\sigma$; in particular, its closure is a link of two 
components. In fact, if we let $\sigma' = \sigma_1^{-2}\sigma$ then 
\[
h(\sigma_1^{-2}\sigma) - h(\sigma) = h(\sigma') - h(\sigma_1^2\sigma') =
-\,(h(\sigma_1^2\sigma') - h(\sigma')),
\]
hence we only need to understand the difference $h(\sigma_1^2\sigma) - 
h(\sigma)$. Let us fix $\ep = (-1, -1, 1, \ldots, 1)$. Since $\sigma_1^2$
and $\ep$ commute, we have 
\begin{multline}\notag
h(\sigma_1^2\sigma) - h(\sigma) 
= \langle\,\w \Delta_n , \w \Gamma_{\ep\sigma_1^2\sigma}\,\rangle 
- \langle\,\w \Delta_n , \w \Gamma_{\ep\sigma}\,\rangle \\
= \langle\,\w \Gamma_{\sigma_1^{-2}} , \w \Gamma_{\ep\sigma}\,\rangle
- \langle\,\w \Delta_n , \w \Gamma_{\ep\sigma}\,\rangle 
= \langle\,\w \Gamma_{\sigma_1^{-2}} - \w \Delta_n , \w \Gamma_{\ep\sigma}\,
\rangle,
\end{multline}
where all intersection numbers are taken in $\w H_n$.
This leads us to consider the {\em difference cycle} $\w\Gamma_{\sigma_1^{-2}} 
- \w \Delta_n$ which is carried by $\w H_n$. The next step in our argument
will be to reduce the analysis of the above intersection to an
intersection theory in the pillowcase $\w H_2$.


\subsection{The pillowcase reduction}
Let us consider the subset $V_n \subset H_n$ consisting of all $(X_1,\ldots, 
X_n, Y_1, \ldots, Y_n) \in H_n$ such that $X_k = Y_k$ for $k = 3,\ldots, 
n$. Equivalently, $V_n$ consists of all $(X_1,\ldots, X_n, Y_1,\ldots, Y_n)
\in Q_n \times Q_n$ such that $(X_1,X_2,Y_1,Y_2) \in H_2$ and $X_k = Y_k$ 
for all $k = 3,\ldots, n$. Therefore, $V_n$ can be identified as
\[
V_n = H_2\times\Delta_{n-2}\;\subset\;(Q_2\times Q_2)\times (Q_{n-2}\times 
Q_{n-2}).
\]

\begin{lemma}
The quotient $\w V_n = (H^*_2 \times \Delta_{n-2})/SO(3)$ is a submanifold 
of $\w H_n$ of dimension $2n - 2$.
\end{lemma}

\begin{proof}
Since $H^*_2$ and $\Delta_{n-2}$ are smooth manifolds of dimensions five and 
$2n-4$, respectively, and their product contains no reducibles, the statement 
follows.
\end{proof}

\begin{lemma}
The difference cycle $\w\Gamma_{\sigma_1^{-2}} - \w \Delta_n$ is carried by 
$\w V_n$.
\end{lemma}

\begin{proof}
Observe that neither $\w \Gamma_{\sigma_1^{-2}}$ nor $\w \Delta_n$ are
subsets of $\w V_n$. However, their portions that do not fit in $\w V_n$, 
\[
\w \Gamma_{\sigma_1^{-2}}
- (\w \Gamma_{\sigma_1^{-2}} \cap \w V_n)\quad\text{and}\quad 
\w \Delta_n - (\w \Delta_n \cap \w V_n),
\]
are exactly the same. Namely, they consist of the equivalence classes of 
$2n$--tuples $(X_1,\ldots,X_n, X_1,\ldots,X_n)$ such that $X_1$ commutes 
with $X_2$. These cancel in the difference cycle $\w\Gamma_{\sigma_1^{-2}} 
- \w \Delta_n$, thus making it belong to $\w V_n$.
\end{proof}

One can isotope $\w\Gamma_{\ep\sigma}$ into $\widetilde\Gamma_{\ep\sigma}$ 
using an isotopy with compact support so that $\widetilde\Gamma_{\ep\sigma}$ 
is transverse to $\w\Gamma_{\sigma_1^{-2}} - \w\Delta_n$. The latter means 
precisely that $\widetilde\Gamma_{\ep\sigma}$ stays away from $(S_2 \times 
\Delta_{n-2})/SO(3)$ and is transverse to both $\w\Gamma_{\sigma_1^{-2}}$ and 
$\w\Delta_n$; a precise argument can be found in \cite[page 491]{HK}. We 
further extend this isotopy to make $\widetilde\Gamma_{\ep\sigma}$ transverse 
to $\w V_n$ so that their intersection is a naturally oriented 1-dimensional 
submanifold of $\w H_n$. 

The natural projection $p: V_n\to H_2$ induces a map $\wp: \w V_n\to \w H_2$.
Use a further small compactly supported isotopy of 
$\widetilde \Gamma_{\ep\sigma}$, if necessary, to make $\wp (\w V_n\,\cap\,
\widetilde \Gamma_{\ep\sigma})$ into a 1-submanifold of $\w H_2$. The proofs 
of Lemmas 2.2 and 2.3 in Lin \cite{Lin} then go through with little change 
to give us the following identity
\[
\langle\, \w \Gamma_{\sigma_1^{-2}} - \w \Delta_n,\; \w \Gamma_{\ep\sigma}\, 
\rangle_{\w H_n} =\, \langle\, \w \Gamma_{\sigma_1^{-2}} - \w \Delta_2,\;
\wp (\w V_n \cap \tilde \Gamma_{\ep\sigma})\, \rangle_{\w H_2}.
\]


\subsection{Computation in the pillowcase}
We begin by studying the behavior of $\wp(\w V_n\,\cap\,\w\Gamma_{\ep\sigma})$ 
near the corners of $\w H_2$. 

\begin{proposition}\label{P:A'}
There is a neighborhood around $A'$ in the pillowcase $\w H_2$ inside which
$\wp (\w V_n\, \cap\, \w \Gamma_{\ep\sigma})$ is a curve approaching $A'$.
\end{proposition}

\begin{proof}
Let us consider the submanifold 
\[
\Delta'_n = \{(X_1,X_2,X_3,\ldots,X_n;Y_1,Y_2,X_3,\ldots,X_n)\}\;\subset\; 
Q_n \times Q_n 
\]
and observe that $V_n\cap \Gamma_{\ep\sigma} = \Delta_n'\cap \Gamma_{\ep\sigma}$.
We will show that the intersection of $\Delta'_n$ with $\Gamma_{\ep\sigma}$ 
is transversal at $(\i,\ep\,\i) = (i,\ldots,i;-i,-i,i,\ldots,i)$. This will 
imply that $\Delta'_n\,\cap\,\Gamma_{\ep\sigma}$ a submanifold of dimension 
four in a neighborhood of $(\i,\ep\,\i)$ and, after factoring out the 
$SO(3)$ symmetry, that $\wp(\w V_n\,\cap\,\w \Gamma_{\ep\sigma})$ is a curve 
approaching $A' = p\,(\i,\ep\,\i)$.

Note that $\dim \Delta'_n = 2n + 4$ hence the dimension of $T_{(\i,\ep\,\i)}
(\Delta'_n\,\cap\,\Gamma_{\ep\sigma}) = T_{(\i,\ep\,\i)}\Delta'_n\,
\cap\,T_{(\i,\ep\,\i)}\Gamma_{\ep\sigma}$ is at least four. Therefore, 
checking the transversality amounts to showing that this dimension is 
exactly four. Write
\[
T_{(\i,\ep\i)}(\Delta_n') = \{(u_1,\ldots u_n;v_1,v_2,u_3,\ldots,u_n)\} 
\subset T_{(\i,\ep\i)}(Q_n \times Q_n)
\]
and
\[
T_{(\i,\ep\i)}(\Gamma_{\ep\sigma}) = \{(u_1,\ldots,u_n; d_{\i}(\ep\sigma)
(u_1,\ldots,u_n)\} \subset T_{(\i,\ep\i)}(Q_n \times Q_n).  
\]
Then $T_{(\i,\ep\i)}(\Delta_n')\,\cap\,T_{(\i,\ep\,\i)}(\Gamma_{\ep\sigma})$ 
consists of the vectors $(u_1,\ldots,u_n) \in T_{\i}\, Q_n = T_i\,S^2\,
\oplus \ldots \oplus\,T_i\,S^2$ that solve the matrix equation 

\begin{equation}\label{E:intsols}
\begin{bmatrix} d_{\i}(\sigma)\end{bmatrix} \;
\begin{bmatrix} u_1 \\ u_2 \\ u_3 \\ \vdots \\ u_n \end{bmatrix}
=
\begin{bmatrix}  * \\  * \\ u_3 \\ \vdots \\ u_n \end{bmatrix};
\end{equation}

\medskip\noindent
since $\ep = (-1,-1,1,\ldots,1)$, we can safely replace $[d_{\i}(\ep\sigma)]$ 
by $[d_{\i}(\sigma)]$. It is shown in \cite{Long} that $[d_{\i}(\sigma)]$ is 
the Burau matrix of $\sigma$ with parameter equal to $-1$. It is a real 
matrix acting on $T_{\i}\, Q_n = \mathbb C^n$, hence all we need to show is 
that the space of $(u_1,\ldots,u_n)\in \mathbb R^n$ solving (\ref{E:intsols}) 
has real dimension two. Let us write 
\[
\begin{bmatrix} d_{\i}(\sigma) \end{bmatrix}
=
\begin{bmatrix} A & B \\ C & D \end{bmatrix} 
\]

\medskip\noindent
where $A$ is a $2 \times 2$ matrix and $D$ is an $(n-2) \times (n-2)$ matrix.
Equation (\ref{E:intsols}) is equivalent to
\[
\begin{bmatrix} C \end{bmatrix}
\begin{bmatrix} u_1 \\ u_2 \end{bmatrix}
=
\begin{bmatrix} 1 - D \end{bmatrix}
\begin{bmatrix} u_3 \\ \vdots \\ u_n \end{bmatrix},
\]
so the proposition will follow as soon as we show that $1 - D$ is invertible.
The invertibility of $1 - D$ is a consequence of the following two lemmas.
\end{proof}

\begin{lemma}\label{E:PermMatrix}
Let $\sigma \in \B_n$ then the Burau matrix of $\sigma$ with parameter $-1$ 
and the permutation matrix of $\bar\sigma^{-1}$ are the same modulo 2. 
\end{lemma}

\begin{proof}
According to \cite{Birman}, the Burau matrix of $\sigma$ with parameter $t$ 
is the matrix 
\[
\frac{\p\sigma(x_i)}{\p x_j}\Big|_{x_i = t}
\]

\medskip\noindent
where $x_i$ are generators of the free group and $\partial$ is the derivative 
in the Fox free differential calculus; see \cite{Fox}. Applying the Fox 
calculus we obtain\,:
\begin{multline}\notag
\frac{\p\sigma(x_i)}{\p x_j}
= \frac{\p(w x_{\bar \sigma^{-1}(i)} w^{-1})}{\p x_j}
= \frac{\p w}{\p x_j} 
+ w \left(\frac{\p(x_{\bar \sigma^{-1}(i)}w^{-1})}{\p x_j}\right) = 
\frac{\p w}{\p x_j} + \\
w \left(\frac{\p x_{\bar \sigma^{-1}(i)}}{\p x_j} 
     + x_{\bar \sigma^{-1}(i)} \frac{\p w^{-1}}{\p x_j}\right)
= \frac{\p w}{\p x_j} 
+ w \frac{\p x_{\bar \sigma^{-1}(i)}}{\p x_j}
- w x_{\bar \sigma^{-1}(i)} w^{-1}\frac{\p w}{\p x_j}, 
\end{multline}

\medskip\noindent
where $w$ is a word in the $x_i$. After evaluating at $t = -1$ and reducing 
modulo 2, the above becomes simply $\p x_{\bar \sigma^{-1}(i)}/\p x_j$, which is 
the permutation matrix of $\bar\sigma^{-1}$.
\end{proof}

\begin{lemma}
Let $\sigma \in \B_n$ be such that $\w\sigma$ is a two component link. Then
$1 - D$ is invertible.
\end{lemma}

\begin{proof}
Our assumption in this section has been that $\bar \sigma = (1 ,\cdots )
(2, \cdots )$. We may further assume that 
\[
\bar \sigma = (1, 3, 4, \ldots, k )\,(2, k+1, k+2, \ldots, n )    
\]
by applying a sequence of first Markov moves fixing the first two strands 
of $\sigma$. The matrix $D \pmod 2$ is obtained by crossing out the first 
two rows and first two columns in the permutation matrix of $\bar\sigma$; 
see Lemma \ref{E:PermMatrix}. This description implies that $D \pmod 2$ is
upper diagonal, and hence so is $(1 - D) \pmod 2$. The diagonal elements 
of the latter matrix are all equal to one, therefore, $\det (1 - D) = 1
\pmod 2$ so $1 - D$ is invertible. 
\end{proof}

\begin{remark}\label{R:orient}
The orientation of the component of $\wp\,(\w V_n\,\cap\,
\w \Gamma_{\ep\sigma})$ limiting to $A'$ can be read off its description near 
$A'$ given in the proof of Proposition \ref{P:A'}. In particular, this 
orientation is independent of the choice of $\sigma$. 
\end{remark}

\begin{proposition}\label{P:AB'}
There are neighborhoods around $A$ and $B'$ in the pillowcase $\w H_2$ which 
are disjoint from $\wp (\w V_n \cap \w \Gamma_{\ep\sigma})$.
\end{proposition}

\begin{proof}
The arguments for $A$ and $B'$ are essentially the same so we will only give 
the proof for $A$. Assuming the contrary we have a curve in $\w V_n \, \cap 
\,\w \Gamma_{\ep\sigma}$ limiting to a reducible representation in $V_n\,\cap
\,\Gamma_{\ep\sigma}$. After conjugating if necessary, this representation 
must have the form 
\[
(i,i, e^{i\phi_3}, \ldots,  e^{i\phi_n}, i, i,  e^{i\phi_3}, \ldots,  e^{i\phi_n}). 
\]
Using the fact that $\ep = (-1, -1, 1, \ldots, 1)$ and arguing as in the proof
of Proposition \ref{P:irrint}, we arrive at the contradiction $i = -i$.
\end{proof}


\subsection{Proof of Theorem 2}
According to Proposition \ref{P:A'}, near $A'$, the 1-submanifold $\wp\,
(\w V_n\, \cap\, \widetilde \Gamma_{\ep\sigma})$ is a curve approaching $A'$. 
According to Proposition \ref{P:AB'}, the other end of this curve must 
approach $B$. Therefore,
\[
h(\sigma_1^2\sigma) - h(\sigma)\, =\, \langle \w \Gamma_{\sigma_1^{-2}} - 
\w \Delta_2,\,\wp (\w V_n\,\cap\,\widetilde \Gamma_{\ep\sigma})\rangle_{\w H_2}
\]
is the same as the intersection number of an arc going from $A'$ to $B$ with 
the difference cycle $\w \Gamma_{\sigma_1^{-2}} - \w \Delta_2$.  This number is 
either $1$ or $-1$ but it is the same for all $\sigma$; see Remark 
\ref{R:orient}. This is sufficient to prove that $h (\sigma)$ is the linking 
number up to an overall sign. 


\end{document}